\documentclass{amsart}

\usepackage{amsmath}
\usepackage{amsxtra}
\usepackage{amscd}
\usepackage{amsthm}
\usepackage{amsfonts}
\usepackage{amssymb}
\usepackage{eucal}

\newtheorem{df}{Definition}[section]

\newtheorem{thm}[df]{Theorem}
\newtheorem{lem}[df]{Lemma}
\newtheorem{claim}[df]{Claim}

\newtheorem{corollary}[df]{Corollary}
\newtheorem{obs}[df]{Observation}

\newtheorem{nota}[df]{Notation}

\def\gA{{\frak A}}

\def\gP{{\frak P}}
\def\gS{{\frak S}}
\def\gK{{\frak K}}
\def\gT{{\frak T}}
\def\gM{{\frak M}}

\def\dcup{{\cup {\!{\!{\!{\cdot}}}}\,}}

\title{Bounded $m$-ary Patch-Width Are Equivalent for $m\geq 3$}
\author{Saharon Shelah and Mor Doron}
\address{The Hebrew university of Jerusalem \newline
            Einstein Institute of Mathematics \newline
            Edmond Safra Campus \newline
            Givat Ram \newline
            Jerusalem 91904, Israel}
\email{shelah@math.huji.ac.il}
\email{mord@math.huji.ac.il}
\thanks{The Author would like to thank the Israel Science Foundation
        for partial support of this research (Grant no. 242/03).}
\thanks{Publication no. 865 in Saharon Shelah's list.}

\begin{document}

\begin{abstract}
We consider the notion of bounded $m$-ary patch-width defined in
\cite{FiMa}, and its very close relative $m$-constructibility
defined below. We show that the notions of $m$-constructibility
all coincide for $m\geq 3$, while $1$-constructibility is a weaker
notion. The same holds for bounded $m$-ary patch-width. The case
$m=2$ is left open.
\end{abstract}

\maketitle

\section{Introduction}
\subsection{Background}
Our interest in this subject started from investigating spectra of
monadic sentences, so let us begin with a short description of
spectra. Let $\phi$ be a sentence in (a fragment of) second order
logic (SOL). The spectrum of $\phi$ is the set $\{n \in
\mathbb{N}:\phi \text{ has a model of size n}\}$. In 1952 Scholz
defined the notion of spectrum and asked for a characterization of
all spectra of first order (FO) sentences. In \cite{As} Asser
asked if the complement of a FO spectrum is itself a FO spectrum.
\begin{df}
    A set $A \subseteq \mathbb{N}$ is eventually periodic if for some
    $n,p \in \mathbb{N}$, for all $m > n$, $m \in A$ iff $m+p \in A$.
\end{df}
In \cite{DuFaLo} Durand, Fagin and Loescher showed that the
spectrum of a FO sentence in a vocabulary with finitely many unary
relation symbols and one function symbol is eventually periodic.
In \cite{GuSh} Gurevich and Shelah generalized this for spectrum
of monadic second order (MSO) sentence in the same vocabulary.
Inspired by \cite{GuSh} Fisher and Makowsky in \cite{FiMa} showed
that the spectrum of a CMSO sentence (a monadic sentence with
counting quantifiers) is eventually periodic provided that all its
models have bounded patch-width. The notion of patch-width of
structures (usually graphs) is a complexity measure on structures,
generalizing clique-width. Their proof remains valid if we
consider $m$-ary patch-width, i.e. we allow $m$-ary relations as
auxiliary relations. In \cite{Sh817} Shelah generalized the proof
of \cite{GuSh} and showed the eventual periodicity of a MSO
sentence provided that all its models are constructible by
recursion using operations that preserve monadic theory (see
definitions below).
\subsection{summation of results}
The above results on eventual periodicity led us to ask: What are
the relations between the the different notions for which we have
eventual periodicity of MSO spectra? In other words do we have
three different results, or are they all equivalent? We give an
answer here. In \cite{Cou} Courcelle proved (using somewhat
different notations) that a class of structures is constructible
iff it is monadicly interpretable in trees, thus implying that two
of the results coincide. We give a proof of Courcelle's result
more coherent with our definition, which we use later on. We prove
that the notions of bounded $m$-ary patch-width is very close to
$m$-constructibility (constructibility where we allow $m$-ary
relations as auxiliary relations) (see lemmas \ref{PW IS CONST}
and \ref{CONST IS PW}). Next we show that for $m\geq 3$ a class of
modes is contained in a $m$-constructible class iff it is
contained in a $3$-constructible class. The same holds for classes
of bounded $m$-ary patch-width (See Theorem \ref{MAIN}). Finally
we show that in the above theorem we can not replace
$3$-constructible by $1$-constructible. That is there exists a
constructible class which is not contained in any
$1$-constructible class. We give a specific example (see
\ref{EXAMPLE}). The case $m=2$ is left open.

\section{Preliminary definitions and previous results}

\begin{nota}  ${}$
\begin{enumerate}

        \item
        Let $\tau$ be a finite relational vocabulary.

        \item
        For $R\in\tau$ let $n(R)$ be the number of places of $R$.
        We say that $R$ is $n(R)$-ary or $n(R)$ place. We allow $n(R)=0$
        i.e. the interpretation of $R$ is in
        $\{\mathbb{T},\mathbb{F}\}$,. We call $\tau$ nice if
        $R\in\tau\Rightarrow n(R)>0$.

        \item
        For $k \in \omega$, let  $\tau_k$ be $\tau \dcup
\{P_{1},...,P_{k}\}$
        with $P_{1},...,P_{k}$ unary predicates.

        \item
        A $k$-colored $\tau$-structure is a $\tau_k$-structure in which
the
        interpretation of the $P_i$'s is a partition of the set of
elements of the
        model (but some $P_i$'s may be empty).

        \item
        A $k$-const $\tau$-structure is a $\tau_k$-structure in
        which every predicate $P_i$ is interpreted by a singleton. We
denote
        such a structure by $(M,a_1,...,a_k)$ where $M$ is a
$\tau$-structure
        and $a_1,...,a_k \in M$.

\end{enumerate}
\end{nota}

\begin{df} ${}$
\begin{enumerate}

    \item
    A monadic second order (MSO) formula in vocabulary $\tau$ is a
    second order formula in which every second order quantifier
    quantifies an \underline{unary} relation symbol. The notion of
    quantifier depth extends, naturally to MSO formulas.

    \item
    Let $M$ be a $\tau$-structure, and $q$ a natural number. The
monadic
    $q$-theory of $M$, $Th_{MSO}^q(M)$, is the set of all sentences
    of quantifier depth $\leq q$ that hold in $M$.

    \item
    Let $M$ be a $\tau$-structure, and $n,q$ natural numbers. Let
    $\bar{a}=(a_1,...,a_n) \in {}^n |M|$. The $q$-type of
    $\bar{a}$ in $M$, $tp_q(\bar{a},M)$, is the set of all
    $\tau$ formulas $\phi$, of quantifier depth $\leq q$ in free
    variables $x_1,...,x_n$, such that: $M \models
    \phi[a_1,...,a_n]$. If $q=0$ we sometimes write
$tp_{qf}(\bar{a}.M)$.

    \item
    The notion of a $q$-type extends to MSO logic. We write
    $tp_q^{MSO}(\bar{a},M)$ for the set of MSO formulas $\phi$,
    of quantifier depth $\leq q$ in free
    variables $x_1,...,x_n$, such that: $M \models
    \phi[a_1,...,a_n]$.

    \item
    The set of all formally possible $q$-types in a vocabulary
    $\tau$ and in variables $\langle x_1,...,x_n \rangle$, will be
    denoted be $TP_q(\langle x_1,...,x_n \rangle,\tau)$, and
    similarly $TP_q^{MSO}(\langle x_1,...,x_n \rangle,\tau)$.
    We may write $TP_q^{MSO}(n,\tau)$ instead of $TP_q^{MSO}(\langle x_1,...,x_n
    \rangle,\tau)$.

\end{enumerate}
\end{df}
\begin{df}[Patch-width]${}$
\begin{enumerate}

\item Let $\tau$ be a nice vocabulary, $M$ a $\tau$-structure, $k$ a natural number, and
$\gP$ a finite set of $k$-colored $\tau$-structures. We say that
$M$ have patch-width at most $k$ (with respect to $\gP$) and
denote $pwd_{\gP}(M) \leq k$, if $M$ is the $\tau$-redact  of a
$k$-colored $\tau$-structure which is in the closer of $\gP$ under
the operations:

\begin{itemize}
        \item[(i)] disjoint union - $\sqcup$,
        \item[(ii)] recoloring - $\rho_{i\to j}$ (change all the elements with color $P_i$ to color $P_j$) and
        \item[(iii)] modifications - $\delta_{R,B}$ (redefine the relation $R\in\tau$ by the quantifier free formula $B$
        in vocabulary $\tau_k$).
\end{itemize}

A class $\gK$ of $\tau$-structures is a $PW(k)$-class, if for some
finite set of $k$-colored $\tau$-structures $\gP$ the elements of
$\gK$ are all the $\tau$-redacts of structures of patch width at
most $k$ with respect to $\gP$. We say $\gK$ is of bounded
patch-width (BPW) if it is a $PW{k}$-class for some
$k\in\mathbb{N}$.

\item In the definition above we may instead of $k$-colored
$\tau$-structures, talk about $\tau^+$ structures where $\tau^+
\supseteq \tau$, $|\tau^+ \setminus \tau| = k$ and every relation
in $\tau^+ \setminus \tau$ is at most $m$-ary. We then talk about
$m$-ary patch-width, where the rest of the definition remains
unchanged. Note that the notions of patch-width and unary patch-
width are close but not identical as in the former we demand that
the colors are disjoint.
\end{enumerate}

\end{df}

In \cite{FiMa} it is proved that:

\begin{thm} \label{PW}
Let $\phi$ be a MSO($\tau$) sentence, and $\gK$ be a $m$-ary
PW(k)-class. Then $spec(\phi)$ restricted to $\gK$ (i.e the set
$\{|M| : M \in \gK , M \models \phi\}$) is eventually periodic.
\end{thm}

\begin{df} [Addition operations] ${}$ \label{AO}
\begin{enumerate}

        \item
        For $k,k_1,k_2 \in \mathbb{N}$, let $\gS_{\tau,k,k_1,k_2}$ be the
set of
        all addition operations of a $k_1$-const $\tau$-structure with
a $k_2$-const
        $\tau$-structure, resulting in a $k$-const $\tau$-structure.
        Formally each $\bf{s}$ $\in \gS_{\tau,k,k_1,k_2}$
        consists of:
        \begin{itemize}

                \item[(i)]
                Sets $A_l=A_l^{\bf{s}} \subseteq \{1,...,k_l\}$ for
                $l \in \{1,2\}$.

                \item[(ii)]
                For $l \in \{1,2\}$, a 1-1 function $g_l=g_l^{\bf{s}}$
from $A_l$ to
                $\{1,...,k\}$ such that: $Im(g_1) \dcup Im(g_2) =
\{1,...,k\}$.

                \item[(iii)]
                For $l \in \{1,2\}$ a set $B_l\subseteq
                \{1,...,k_l\}^2$, and a set $B \subseteq
                \{1,...,k_1\}\times\{1,...,k_2\}$.

                \item[(iv)]
                For each $R \in \tau$ with $n(R)=n$
                and each $w_l \subseteq \{1,...,n\}$ for $l \in
\{1,2\}$,
                a function $f_{R,w_1,w_2}=f_{R,w_1,w_2}^{\bf{s}}$
                with range $\{\mathbb{T},\mathbb{F}\}$, and domain: triplets of
the form
                $(p,q_1,q_2)$ where:
                \begin{itemize}

                        \item
                        $p\in TP_0(\langle x_1,...,x_n \rangle,\sigma)$ were $\sigma$ is
                        a vocabulary with $k_1+k_2$ individual constants and two unary predicates,

                        \item
                        For $l \in \{1,2\}$, $q_l \in TP_0(\langle x_i
: i \in w_l \rangle ,\tau)$.

                \end{itemize}

        \end{itemize}

        \item
        Let  $k,k_1,k_2 \in \mathbb{N}$ and  $\bf{s}$ $\in
\gS_{\tau,k,k_1,k_2}$. Let
        $(M_l,a^l_1,...,a^l_{k_l})$ be $k_l$-const $\tau$-structure for
$l \in \{1,2\}$.
        The addition
        $(M_1,a^1_1,...,a^1_{k_1}) \circledast_{\bf{s}}
(M_2,a^2_1,...,a^2_{k_2})$
        is defined whenever:
        \begin{itemize}
        \item
        $(|M_1| \cap |M_2|) \subseteq (\{a^1_1,...,a^1_{k_1}\} \cap
        \{a^2_1,...,a^2_{k_2}\})$ and
        \item
        For $l \in \{1,2\}$: $a^l_i=a^l_j\Leftrightarrow (i,j)\in
        B_l$ and $a^1_i=a^2_j\Leftrightarrow (i,j)\in
        B$,
        \end{itemize}
        to be the $k$-const $\tau$-structure $(M,b_1,...,b_k)$ defined
by:
        \begin{itemize}

                \item[(i)]
                $|M|=(|M_1| \setminus \{a^1_1,...,a^1_{k_1}\}) \cup
(|M_2| \setminus \{a^2_1,...,a^2_{k_2}\})
                 \cup \{a_i^l : l \in \{1,2\}, i \in A_l\})$.

                \item[(ii)]
                For each $l \in \{1,2\}$ and $i \in A_l$,
$a_i^l=b_{g_l(i)}$.

                \item[(iii)]
                For all $R \in \tau$ with $n(R)=n$ and
                $\bar{x}=(x_1,...,x_n) \in {}^n|M|$, let
                $w_l=\{ i:x_i \in |M_l| \}$ for $l \in \{1,2\}$.
                Let $p$ be the quantifier free type of $\bar{x}$ in the
model
                with $\{a^1_1,...,a^1_{k_1}\} \cup
\{a^2_1,...,a^2_{k_2}\}$ as constants,
                and $|M_1|, |M_2|$ as unary predicates.
                For $l \in \{1,2\}$ let $q_l=tp_{qf}(\langle x_i: i \in
w_l \rangle,M_l)$.
                Now the value of $R^M(\bar{x})$
                is defined to be $f_{R,w_1,w_2}^{\bf{s}}(p,q_1,q_2)$.

        \end{itemize}

        \item
        For technical reasons we would like to allow empty structures.
        i.e.  let $\tau':=\{R\in\tau:n(R)=0\}$, and $X\subseteq\tau'$. Now
        $Null_X$ is the $\tau'$-structure with $|Null_X|=\emptyset$ and
        $R^{Null_X}=True\Leftrightarrow R\in X$.
        Then if $\bf{s}$ $\in \gS_{\tau,k,k_1,0}$,
        and $M$ is a $\tau_{k_1}$-
        structure then $M \circledast_{\bf{s}} Null_X$ is a well defined
        $\tau_k$-structure. Furthermore for any $\tau$-structure
        $M$, $M\sqcup Null_{\emptyset}$ is defined and equal to $M$.

\end{enumerate}
\end{df}

The important attributes of the addition operations are the
following:

\begin{thm}
Let $k,k_1,k_2 \in \mathbb{N}$. Then:
\begin{enumerate}

        \item $\gS_{\tau,k,k_1,k_2}$ is finite.

        \item
        \emph{The addition theorem:} \newline
        Let $M,M'$ be $k_1$-const $\tau$-structures such that
$Th^q_{MSO}(M)=Th^q_{MSO}(M')$,
        and $N,N'$ be $k_2$-const $\tau$-structures such that
$Th^q_{MSO}(N)=Th^q_{MSO}(N')$,
        and $\bf{s}$ $\in \gS_{\tau,k,k_1,k_2}$. Assume that the
additions
        $M \circledast_{\bf{s}} N$ and  $M' \circledast_{\bf{s}} N'$
are defined. Then
        \[Th^q_{MSO}(M \circledast_{\bf{s}} N) = Th^q_{MSO}(M'
\circledast_{\bf{s}} N').\]

\end{enumerate}
\end{thm}

\begin{df} [Constructibility]\label{def const}
A class $\gK$ of $\tau$-structures is $(m^*,k^*)$-constructible,
if there exists: A finite relational vocabulary $\tau^+ \supseteq
\tau$, a finite set of structures $\gP$, and a finite set of
addition operations $\gS$ such that:

\begin{itemize}
        \item[(i)]
        Every relation in $\tau^+ \setminus \tau$ is at most $m^*$-ary.
        \item[(ii)]
        Every structure in $\gP$ is a $k$-const $\tau^+$-structure
        for some $k \leq k^*$.

        \item[(iii)]
        Every operation in $\gS$ is in $\gS_{\tau^+,k,k_1,k_2}$ for
some
        $k,k_1,k_2 \leq k^*$.

        \item[(iv)]
        The elements of $\gK$ are all the $\tau$-redacts of structures
in
        the closer of $\gP$ under the operations in $\gS$.
\end{itemize}
We say that $\gK$ is $m^*$-constructible if it is $(m^*,k^*)$-
constructible for some $k^*$.

\end{df}

In \cite{Sh817} we see that:

\begin{thm}
Let $\phi$ be a MSO($\tau$) sentence, and $\gK$ be a
$(m^*k^*)$-constructible class. Then $spec(\phi)$ restricted to
$\gK$ (i.e the set $\{|M| : M \in \gK , M \models \phi\}$) is
eventually periodic.
\end{thm}

This is a generalization of \ref{PW} as we have:

\begin{lem} \label{PW IS CONST}
Let $\tau$ be a nice vocabulary, and $\gK$ be a $m$-ary
$PW(k)$-class of $\tau$-structures. Then $\gK$ is a
$(m,0)$-constructible class.
\end{lem}

\begin{proof}
First note that the disjoint union operation of
$\tau^+$-structures is in $\gS_{\tau^+,0,0,0}$. As for the
recoloring and the modification operations, those are unary
operations, so we look at the operation ${\bf{s}} \in
\gS_{\tau^+,0,0,0}$ that acts as recoloring or modification on its
left operand. So $M \circledast_{\bf{s}} Null_{\emptyset}$ is the
desired recoloring or modification of $M$.
\end{proof}

In the addition operations we allow omitting marked elements, and
the universe of the the two operands is not necessarily disjoint.
This is not allowed in the operations of patch-width. It turns out
though that these are the only essential differences between the
two types of operations as suggested by the following:

\begin{lem} \label{CONST IS PW}
Let $\gK$ be a $(m,0)$-constructible class such that the set of
operations $\gS$ associated with $\gK$ satisfies
$\gS\subseteq\gS_{\tau^+,0,0,0}$, and the vocabulary $\tau^+$
associated with $\gK$ is nice. Then $\gK$ is am $m$-ary $PW(k')$
class for $k'=2(k+|\tau|)+2$.
\end{lem}

\begin{proof}
$\gK$ is $(m,k)$-constructible so we have a vocabulary $\tau^+$
and sets $\gS$ and $\gP$. Now the set of atomic structures for the
patch-width definition will be the same $\gP$. The vocabulary of
the patch-width definition will be: $\tau^+ \dcup
\{R':R\in\tau^+\} \dcup \{P_1,P_2\}$. $P_1,P_2$ are new unary
relation symbols, note that the number of relations in the new
vocabulary is indeed $k'$. We now have to show for each operation
in $\gS$ how to simulate it by operations of patch-width. Let
$\bf{s}\in\gS$ and let $M_1,M_2$ be $\tau^+$-structures. Denote by
$M_1',M_2'$ the trivial extensions to the new vocabulary. We will
now describe a series of patch-width operations on $M_1',M_2'$
resulting in a structure $M^*$ such that $M^*|_{\tau^+} \cong M_1
\circledast_{\bf{s}} M_2$, this will complete the proof. First
color all the elements of $M_l'$ by $P_l$ for $l\in\{1,2\}$. Next
for each $R\in\tau^+$ redefine $R'$ to be the same as $R$, do this
for both $M_1',M_2'$. Now take the disjoint union of the to
resulting structures. Finally we have to redefine the relations of
$\tau^+$ of the our disjoint union to be as in $M_1
\circledast_{\bf{s}} M_2$. Let $R\in\tau^+$ be $n$-ary and let
$w_1,w_2\subseteq \{1,...,n\}$ satisfy $w_1\dcup w_2=\{1,...,n\}$.
Let $p$ be the quantifier free type in the vocabulary with two
unary relations $S_1,S_2$ "saying" that for $i\leq n$ and
$l\in\{1,2\}$, $x_i\in S_l$ iff $i\in w_l$. Now define:
$$\varphi_{R,w_1,w_2}(x_1,...x_n):=\bigwedge_{i\in
w_1}P_1(x_i)\bigwedge_{i\in w_2}P_2(x_i) \wedge
  [\bigvee_{\substack {q_l\in TP_0(\langle x_i:i\in w_l\rangle ,\tau^+)
\\ f^{\bf{s}}_{R,w_1,w_2}(p,q_1,q_2)=\mathbb{T}}}
   \wedge q_1'\wedge q_2'].$$
Where $\wedge q_l'$ is the disjunction of all the formulas in
$q_l$ where we replace every relation $R\in \tau^+$ be $R'$. Now
redefine the relation $R$ using the modification $\delta_{R,B}$
for the formula:
$$B(x_1,...,x_n):=\bigwedge_{\substack{w_1,w_2\subseteq \{1,...,n\} \\
w_1\dcup w_2=\{1,...,n\}}}
  \varphi_{R,w_1,w_2}(x_1,...,x_n).$$
Do this for all $R\in\tau^+$ and we are done.
\end{proof}


\begin{nota}[Trees]  ${}$
\begin{enumerate}

        \item
        The vocabulary of trees, $\tau_{trees}$, is $\{\leq,c_{rt}\}$.

        \item
        The vocabulary of $k$-trees, $\tau_{k-trees}$, is
        $\{\leq,c_{rt}\} \dcup \{P_1,...,P_k\}$ i.e
${(\tau_{trees})}_k$.

        \item
        A tree $\gT$ is a $\tau_{trees}$-structure in which:

        \begin{itemize}
                \item
                For every $t \in |\gT|$ the set $\{s \in |\gT| : s
\leq^{\gT} t\}$
                is linearly ordered by $\leq^{\gT}$.
                \item
                For all $x\in |\gT|$, ${c_{rt}}^{\gT}\leq^{\gT}x$.
        \end{itemize}

        \item
        A $k$-tree $\gT$ is a $\tau_{k-trees}$-structure,
        such that $\gT|\tau_{trees}$ is a tree.

        \item
        A $2$-tree $\gT$ is directed binary (DB) if
$(c_{rt}^{\gT},P_1^{\gT},P_2^{\gT})$
        is a partition of $|\gT|$, and each non maximal element of
$\gT$ has exactly two
        immediate successors one in $P_1^{\gT}$ and the other in
$P_1^{\gT}$. For $k\geq 2$,
        a $k$-tree $\gT$ is DB if
$(|\gT|;\leq^{\gT},c_{rt}^{\gT},P_1^{\gT},P_2^{\gT})$ is.

\end{enumerate}
\end{nota}

\begin{df}[Monadic interpretation] ${}$
\begin{enumerate}
        \item
        We call $\bf{c}$ a  monadic $k$-interpretation scheme for a
vocabulary $\tau$
        if $\bf{c}$ consists of:

        \begin{itemize}
                \item
                Natural numbers $k_1=k_1^{\bf{c}}$ and
$k_1=k_2^{\bf{c}}$ both less then or equal to $k$.

                \item
                For every $l \leq k_1$ a monadic
$\tau_{k_2-trees}$-formula $\varphi_{=,l}^{\bf{c}}(x)$.

                \item
                For every $R \in \tau$ $n$-place relation, and every
                $\eta \in {}^{\{1,...,n\}}\{0,...,k_1\}$ a monadic
$\tau_{k_2-trees}$-formula:
                $\varphi=\varphi^{\bf{c}}_{R,\eta}(x_1,...,x_n)$.
        \end{itemize}

        \item
        Let $\bf{c}$ be a monadic $k$-interpretation scheme for a
vocabulary $\tau$,
        and $\gT$ a $k_2^{\bf{c}}$-tree. The interpretation of $\gT$ by
        $\bf{c}$ denoted by $\gT^{[\bf{c}]}$ is the $\tau$-model $M$
defined by:
        \begin{itemize}
                \item
                $|M|=\{(t,l) \in |\gT| \times \{0,...,k_1\} : \gT
\models \varphi_{=,l}(t)\}$

                \item
                For every $R \in \tau$ $n$-place relation:
                $$R^M=\{((t_i,l_i):i \leq n) \in {}^n |M|: \gT \models
                \varphi_{R,(l_i:i \leq n)}(t_1,...,t_n)\}$$
        \end{itemize}

        \item
        A $\tau$-model $M$ is monadicly $k$-interpretable in trees if
for some $\bf{c}$ a
        monadic $k$-interpretation scheme for $\tau$, and some
$k_2^{\bf{c}}$-tree,  $\gT$,
        we have: $\gT^{[\bf{c}]} \cong M$. We denote the class of all the
        $\tau$-structures
        monadicly $k$-interpretable in trees by $\gK^{mo}_{\tau,k}$.

        \item
        For $\bf{c}$ a monadic $k$-interpretation scheme for $\tau$ we
        denote by
        $\gK_{{\bf{c}}}^{mo}$ the class of all $\tau$-structures
        $M$ such that for some $k_2^{\bf{c}}$-tree,  $\gT$,
        we have: $\gT^{[\bf{c}]} \cong M$. $\gK_{{\bf{c}}}^{mo,db}$ is
        the same as
        $\gK_{{\bf{c}}}^{mo}$ only we demand that $\gT$ is directed
        binary.

        \item
        We say that $\bf{c}$ has the leaf property if
        $k_1^{\bf{c}}=0$ and for every $k_2^{\bf{c}}$-tree $\gT$,
        and every $t\in|\gT|$, we have:
        $\gT\models\varphi_{0,=}^{\bf{c}}[t]$ implies that $t$ is
        a maximal element in $\gT$.

\end{enumerate}
\end{df}

Without loss of generality we may assume that $k_1^{\bf{c}}=0$.
This is because of the following:

\begin{lem} \label{wlog}
For every $\bf{c}$ a monadic $k$-interpretation scheme for a
vocabulary $\tau$, there exists $\bf{c'}$ be a monadic
$(k+2)$-interpretation scheme for $\tau$, such that:
    \begin{itemize}
    \item $k_1^{\bf{c'}}=0$.
    \item $k_2^{\bf{c'}}=k_2^{\bf{c}}+2$.
    \item For every $k_2^{\bf{c}}$-tree, $\gT$ there exists a
$k_2^{\bf{c'}}$-tree,
    $\gT'$, such that: $\gT^{[\bf{c}]} \cong \gT'^{[\bf{c'}]}$.
    \end{itemize}
Hence  $\gK_{{\bf{c}}}^{mo} \subseteq \gK_{{\bf{c'}}}^{mo}$.
\end{lem}

\begin{proof}
Let $s_1$ and $s_2$ be the two "new" unary predicates, and let
$\gT$ be a $k_2^{\bf{c}}$-tree. Define $\gT'$ as follows:
$|\gT'|=|\gT| \dcup (|\gT| \times \{0,...,k_1^{\bf{c}}\})$,
$s_1^{\gT'}=|\gT|$, $s_2^{\gT'}= |\gT| \times
\{0,...,k_1^{\bf{c}}\}$, and if $t_1$ is the immediate successor
of $t_2$ in $\gT$ then define, $t_1 <^{\gT'} (t_1,0) <^{\gT'}
(t_1,1) <^{\gT'} ... (t_1,k_1^{\bf{c}}) <^{\gT'} t_2$. Now define:
    $$\varphi_{=,0}^{\bf{c'}}(x):=s_2(x) \wedge
    \bigwedge_{l<k_1^{\bf{c}}} (\forall y) [s_1(y) \wedge (\psi_l(x,y)]
\to
    (\varphi_{=,l}^{\bf{c}}(y))^{s_1}$$

Where $\psi_l(x,y)$ is a formula stating that there are exactly
$l$ elements between $x$ and $y$ and all of them are in $s_2$, and
$(\varphi_{=,l}^{\bf{c}}(y))^{s_1}$ is the formula
$\varphi_{=,l}^{\bf{c}}(y)$ relativized to $s_1$ i.e we replace
every quantifier of the form $\exists x$ or $\forall x$ by
$\exists x \in s_1 $ or $\forall x \in s_1$ respectively. It
should be clear that $\gT \models \varphi_{=,l}^{\bf{c}}[t]$ iff
$\gT' \models \varphi_{=,0}^{\bf{c'}}[(t,l)]$. The relations are
dealt with in a similar way.
\end{proof}

\begin{lem} \label{CONST IS TREE}
Let $\gK$ be a $(m^*,k^*)$-constructible class of $\tau$-models.
Then there exists a natural number $k^{**}$ such that $\gK
\subseteq \gK^{mo}_{\tau,k^{**}}$. Moreover for some monadic
$k^{**}$-interpretation scheme $\bf{c}$ with the leaf property, we
have $\gK \subseteq \gK^{mo,db}_{{\bf{c}}}$.
\end{lem}

We will not go into detail here especially as a similar result was
proved by Courcelle in \cite{Cou}. We do however give a sketch of
a proof containing some definitions that will be useful later.

\begin{proof}[Sketch]
Suppose $\gP$ and $\gS$ are the finite sets of structures and
operations generating $\gK$, and $\tau^+$ the vocabulary
associated with $\gK$ (see \ref{def const}). Now with every $M \in
\gK$ we can associate a tree which represents the construction of
$M$ from the structures in $\gP$. Formally we define:

        \begin{df} \label{rep}
        We say that the pair $(\gT,\gM)$ with
        $\gT=\langle T;\leq^{\gT},c_{rt}^{\gT},S_1^{\gT},S_2^{\gT}
\rangle$
        a DB tree, and
        $\gM=\langle M_t : t \in T \rangle$
        is a full representation of $M \in \gK$ when:

        \begin{enumerate}
                \item
                Every $M_t$ is a $k_t$-const $\tau^+$-structure for
some $k_t \leq k^*$.

                \item
                For every $t \in T$ $\leq^{\gT}$-maximal, $M_t \in
\gP$.

                \item
                The $\tau$-redact of $M_{c_{rt}^{\gT}}$ is $M$.

                \item
                For every $t$, a non-maximal element of $T$, let
$s_1,s_2$
                be its immediate successors with $s_l \in S_l^{\gT}$.
Then
                $M_t=M_{s_1} \circledast_{\bf{s}} M_{s_2}$ for some
                ${\bf{s}} \in \gS_{\tau^+
                ,k_{s_1},k_{s_2},k_t} \cap \gS$.

        \end{enumerate}

        \end{df}

        \begin{df} ${}$ \label{def:tau*}
        \begin{enumerate}

                \item
                Let $\tau^*$ be the vocabulary $\tau_{k_2-trees}$ with
the following
                unary predicates:
                \begin{enumerate}
                        \item $S_1$ and $S_2$.
                        \item $P_k$ for $k \leq k^*$.
                        \item $Q_{\bf{s}}$ for ${\bf{s}} \in \gS$.
                        \item $R_N$ for $N \in \gP$.
                \end{enumerate}
                $k_2$ is the total number of unary predicates in
                $\tau^*$, i.e $k_2=|\gP|+|\gS|+k^*+2$.

                \item
                A $\tau^*$-structure $\gT$ is a representation of $M
\in \gK$,
                if we can find         $\gM=\langle M_t : t \in |\gT|
\rangle$
                such that:
                \begin{enumerate}
                        \item

$((|\gT|,\leq^{\gT},c_{rt}^{\gT},S_1^{\gT},S_2^{\gT}),\gM)$
                        is a full representation of $M$.

                        \item
                        $\langle P_k^{\gT} : k \leq k^* \rangle$ is a
partition of $|\gT|$.
                        If $t \in P_k^{\gT}$, then $k_t = k$ i.e. $M_t$
is a $k$-const $\tau^+$-structure.
                        We write $k^{\gT}(t)=k$ iff $t \in P_k^{\gT}$.

                        \item
                        $\langle Q_{\bf{s}}^{\gT} : {\bf{s}} \in \gS
\rangle \cup
                        \langle R_N^{\gT} : N \in \gP \rangle$ is a
partition of $|\gT|$.

                        \item
                        For every $t \in T$ $\leq^{\gT}$-maximal, $t
\in R_{M_t}^{\gT}$.

                        \item
                        For every $t \in T$ non-maximal, let $s_1,s_2$
                        be its immediate successors with $s_l \in
S_l^{\gT}$. Suppose
                        $M_t=M_{s_1} \circledast_{\bf{s}} M_{s_2}$ for
some
                        ${\bf{s}} \in \gS_{\tau^+,k_{s_1},k_{s_2},k_t}
\cap \gS$. Then
                        $t \in Q_{\bf{s}}^{\gT}$.

                \end{enumerate}

        \end{enumerate}
        \end{df}

        Note that:

        \begin{obs} ${}$
        \begin{enumerate}
                \item
                Every $M \in \gK$ has a full representation, and
                hence a representation.

                \item
                If $M_l \in \gK$ are represented by $\gT_l$ for
$l\in\{1,2\}$,
                and $\gT_1 \cong \gT_2$. Then $M_1 \cong M_2$.
        \end{enumerate}
        \end{obs}

        Now define: $k_1=max\{|N|:N \in \gP\}$,
        $k_2$ is the number of unary predicates in $\tau^*$ (see
        \ref{def:tau*}(1)),
        and let $k^{**}=max\{k_1,k_2\}$. We can define a
        $k^{**}$-interpretation scheme
        $\bf{c}$ with $k_1^{\bf{c}}=k_1$ and  $k_2^{\bf{c}}=k_2$ such
        that for all
        $M \in \gK$, and $\gT$ a representation of $M$ we have  $M
        \cong \gT^{[\bf{c}]}$.
        Note that indeed $\gT$ is a DB $k^{\bf{c}}_2$-tree. We
        will not specify all the formulas of $\bf{c}$ as they tend to be very long and
        complicated, but do note that all the information about $M$
        can be decoded from the representation of $M$ using
        monadic formulas. Finally by an argument very close to
        that of \ref{wlog} we may assume that $\bf{c}$ has the
        leaf property.
\end{proof}


\section{Equivalence of $m$-ary patch-width For $m\geq 3$}

We come now to the main part of our result. Basically what we do
here is proving the reverse inclusion of \ref{CONST IS TREE}. It
turns out that in our constructible class we only need $3$-ary
relations as auxiliary relations, thus we can replace
constructible by $(3,k)$-constructible. It follows that a class
$\gK$ is contained in a $(m,k)$-constructible class for some $k$,
iff it is contained in a $(3,k)$-constructible class for some $k$,
and similarly for $m$-ary path-width. We start with an
investigation of directed binary trees that will be useful later.

\begin{nota}\label{nota1}
Let $\gT$ be a DB $k$-tree. Let $n\in\mathbb{N}$ and
$x_1,...,x_n\in T$ be fixed maximal elements of $\gT$.
\begin{enumerate}
\item For $x,y\in T$ denote by $x\wedge y$ the minimal element $z$
with $z\geq x,y$.
\item
For $x,y\in T$ with $x\leq y$ denote
$[x,y):=\{z\in T x\leq z<y\}$ and similarly
      $(x,y),(x,y]$ and $[x,y]$.
\item Define $Y:=\{x_1,...,x_n\}\cup\{x_i\wedge x_j:i,j\leq
n\}\cup\{c_{rt}^{\gT}\}$. Note that $|Y|\leq 2n$. \item For any
non-maximal $x\in T$ let $F_R(x)\in T$ (resp. $F_L(x)\in T$) be
the unique immediate successor of $x$ which is in $P_1^{\gT}$
(resp. $P_2^{\gT}$). \item For $y,y'\in Y$ with $y<y'$ define,
$$\begin{array}{llr}
                                              T^3_{y,y'}:= & [y,y']\cup \\
                                               & \{z\in T:(\exists s\in
(y,y'))F_R(s)\leq y' \wedge F_L(s)\leq z\} \cup \\
                                               & \{z\in T:(\exists s\in
(y,y'))F_L(s)\leq y' \wedge F_R(s)\leq z\}
\end{array}$$
\item Let $T_R=\{c_{rt}^{\gT}\}\cup\{t\in T:F_R(c_{rt}^{\gT})\leq
t\}$, and similarly $T_L$.
\end{enumerate}
\end{nota}

\begin{lem} \label{const}
Let $R_R(y,y')$ and $R_L(y,y')$ be binary relations meaning
$F_R(y)\geq y'$ and $F_L(y)\geq y'$ respectively. The type
$tp_q^{MSO}((x_1,...,x_n),\gT)$ is computable from the structure
$\langle Y; \leq^{\gT}, R_R, R_L \rangle$, the types
$\{tp_q^{MSO}((y,y'),\gT|_{T^3_{y,y'}}):y,y'\in Y,y<y',(y,y') \cap
Y = \emptyset\}$, and the types
$tp_q^{MSO}(c_{rt}^{\gT},\gT|_{T_L})$,
$tp_q^{MSO}(c_{rt}^{\gT},\gT|_{T_R})$.
\end{lem}

\begin{proof}
Without going into detail note that from the sets $T^3_{y,y'}$
with $y,y'$ as above, $T_L$ and $T_R$, we can choose a
decomposition of $|\gT|$, in which only the elements of $Y$ belong
to more then one set. Hence we can reconstruct the structure $\gT$
with the elements of $Y$ as marked elements from the reduced
structures: $\gT|_{T^3_{y,y'}}$ with $y,y'$ as marked elements,
$\gT|_{T_L}$ and $\gT|_{T_L}$ with $C_{rt}^{\gT}$ as marked
element, in a way that the $q$ theory of the resulting structure
depends only on the $q$ theory of the operands. The structure
$\langle Y; \leq^{\gT}, R_R, R_L \rangle$ determines the order of
the construction.
\end{proof}

\begin{claim} \label{tree const}
        Let $k^*$ be a natural number, and $\bf{c}$ a
        monadic $k^*$-interpretation scheme with the leaf property
        for a vocabulary $\tau$.
        Then there exists a natural number $k^{**}$, and a
$(3,k^{**})$-constructible
        class of $\tau$ structures, $\gK$, such that:
        $\gK^{mo,db}_{{\bf{c}}} \subseteq \gK$.
\end{claim}

\begin{proof}
    Let $q^*$ be the maximal quantifier rank of the formulas
$\{\varphi_{Q,0}:Q \in\tau\}$.
    Define the vocabulary $\tau^+$ to consist of:
    \begin{itemize}
        \item $\tau$.
        \item $\tau_{k_2-trees}$.
        \item Two $3$-place relations $R_R$ and $R_L$.
        \item
            For each $\textbf{t} \in
TP_{q^*}^{MSO}(2,\tau_{k_2-trees})$,
            a $3$-place relation $R^3_{\textbf{t}}$.
            \item
            For each $\textbf{t} \in
TP_{q^*}^{MSO}(2,\tau_{k_2-trees})$,
            a $2$-place relation $R^2_{\textbf{t}}$.
        \item For each $\textbf{t} \in
TP_{q^*}^{MSO}(1,\tau_{k_2-trees})$ two $0$-place relations
$R^R_{\bf{t}}$ and $R^L_{\bf{t}}$.

    \end{itemize}

    Before we define the set of addition operations $\gS$, and the set
$\gP$, let us
    define:

    \begin{df} ${}$ \label{correct}
        A $\tau^+$-structure, $\gT$, is called a "correct" $k_2$-tree
if:
        \begin{itemize}
            \item For each $Q \in \tau$, $Q^{\gT}=\emptyset$.
            \item $\gT|\tau_{k_2-trees}$ is a DB $k_2$-tree.
            \item
                For each $x_1,x_2,x_3$ maximal elements of $|\gT|$, let
                $y=x_1\wedge x_2$ and $y'=y\wedge x_3$ then we have,
                $R_{R}^{\gT}(x_1,x_2,x_3)\Leftrightarrow F_R(y)\geq
y'$,
                and similarly for $R_L$.
                 \item
                For each $\textbf{t} \in
TP_{q^*}^{MSO}(2,\tau_{k_2-trees})$,
                and $x_1,x_2,x_3$ maximal elements of $|\gT|$, let
                $y=x_1\wedge x_2$ and $y'=y\wedge x_3$ then we have,
                $(R^3_{\textbf{t}})^{\gT}(x_1,x_2,x_3)\Leftrightarrow
                tp_{q^*}^{MSO}((y,y'),\gT|_{T^3_{y,y'}}) = \textbf{t}$.
                 \item
                For each $\textbf{t} \in
TP_{q^*}^{MSO}(2,\tau_{k_2-trees})$,
                and $x_1,x_2$ maximal elements of $|\gT|$, let
                $y=x_1\wedge x_2$ then we have,
                $(R^2_{\textbf{t}})^{\gT}(x_1,x_2)\Leftrightarrow
                tp_{q^*}^{MSO}((c_{rt}^{\gT},y),\gT|_{T^3_{c_{rt}^{\gT},y}}) = \textbf{t}$.
            \item
                For each $\textbf{t} \in
TP_{q^*}^{MSO}(1,\tau_{k_2-trees})$,
$(R^R_{\bf{t}})^{\gT}=\mathbb{T}$ iff
$tp_{q^*}^{MSO}(c_{rt}^{\gT},\gT|_{T_R})=\bf{t}$, and similarly
for $R^L_{\bf{t}}$.
        \end{itemize}
    \end{df}

    Note that every DB $k_2$-tree can be uniquely extended to a
    correct DB $k_2$-tree.
    Now define Our $\gP$ to consist of all singleton correct models
    (models with one element) of the vocabulary $\tau^+$, plus all the Null
    $\tau^+$-structures (see definition \ref{AO}(5)) .

    We now turn to the definition of the operations in $\gS$.
    Let $u$ be a possible "color" of a singleton $k_2$ tree. Formally
    $u\subseteq\{P_3,...,p_{k_2}\}$. We define the operation $\oplus_u$
    on DB $k_2$-trees as the addition of two trees with root of color $u$. Formally
    Let $\gT_1, \gT_2$ be DB $k_2$-trees define $\gT=\gT_1\oplus_u\gT_2$ by:
    \begin{itemize}
        \item $|\gT|=|\gT_1| \dcup |\gT_2| \dcup \{c\}$.
        \item $c$ is the root of $\gT$ i.e. $c_{rt}^{\gT}=c$ and $\forall t\in |\gT| c<^{\gT} x$.
        \item $c$ has color $u$ i.e. for all $i\geq 3$, $c\in P_i^{\gT}$ iff $i\in u$.
        \item $c_{rt}^{\gT_1}\in P_1^{\gT}$ and $c_{rt}^{\gT_2}\in
        P_2^{\gT}$.
        \item The rest of the relations on $\gT_1$ and $\gT_2$
        remain unchanged.
    \end{itemize}
    Note that indeed $\gT_1\oplus_u\gT_2$ is a DB $k_2$-tree whenever $\gT_1$ and
    $\gT_2$ are, and hence $\oplus_u$ extends uniquely to an operation on
    \underline{correct} $k_2$-trees.

    Now for $l\in\{1,2\}$ let $\gA_l$ be a $\tau^+$ structure such
    that there exists a correct $k_2$-tree with $|\gA_l|\subseteq
    |\gT_l|$, $\gT_l|_{|\gA_l|}=\gA_l$, and every element of
    $\gA_l$ is maximal in $\gT_l$. Define am operation $\bf{s}_u$ on such structures
    by: $\gA_1 \circledast _{{\bf{s}}_u} \gA_2 = (\gT_1 \oplus_u
    \gT_2)|_{|\gA_1|\cup |\gA_2|}$. It is easy to verify that
    $\circledast_{{\bf{s}}_u}$ is well defined and indeed belongs
    to $\gS_{\tau^+,0,0,0}$. We now have:

    \begin{lem} \label{ind}
        For every correct $k_2$-tree, $\gT$ and every set $A\subseteq |\gT|$
        of maximal elements, the restriction $\gT|_A$
        is in the closer of $\gP$
        under the operations $\{\bf{s}_u:u\subseteq\{1,...,k_2\}\}$.
    \end{lem}

    \begin{proof}
        First it is obvious that we can construct $\gT$ from $\gP$
        using the operations $\{\oplus_u:u\subseteq\{1,...,k_2\}\}$.
        Now use the same construction only replace the operation
        $\gT_1\oplus_u\gT_2$ by the operation
        $\gT_1|_A\circledast_u\gT_2|_A$.
    \end{proof}

    The last thing we need now is to "decode" the relations
    in the correct structure into the relations in
    our vocabulary $\tau$. For this we use:

    \begin{lem} \label{s*}
    There exist ${\bf{s}^*} \in \gS_{\tau^+,0,0,0}$ such that For
    every correct $k_2$-tree, $\gT$ and every set $A\subseteq |\gT|$
    of maximal elements, the structure $\gA'=\gT|_A\circledast_{\bf{s}^*}NULL_{\emptyset}$
    satisfies for each $Q\in \tau$ with $n(Q)=n$,

        $$(*)\quad Q^{\gA'} = \{(x_1,...,x_n) \in {}^n A : \gT \models
        \varphi_{Q,0}(x_1,...,x_n)\}.$$

    \end{lem}

    \begin{proof}
    Let $Q \in \tau$ be an $n$-place relation symbol, and $w_1,w_2
    \subseteq \{1,...,n\}$.
    We should define $f_{Q,w_1,w_2}^{\bf{s_3}}$ in such a way that
    $(*)$ will hold. As we have $k_1^{{\bf{s}}^*}=k_2^{{\bf{s}}^*}=k^{{\bf{s}}^*}=0$
    and we are only interested in $NULL_{\emptyset}$ as the left operand,
    the only relevant
    case is $w_1=\{1,...,n\}$ and $w_2=\emptyset$.
    In order to have $(*)$ We need to define a function:
    $$ f_Q^{{\bf{s}^*}} : \{p: \text{p is a quantifier free type of n variables in
vocabulary } \tau^+\}
    \rightarrow \{\mathbb{T},\mathbb{F}\}$$
    such that for all $(x_1,...,x_n) \in {}^nA$,
$f_Q(tp_{qf}((x_1,...,x_n),\gA'))=\mathbb{T}$ iff
    $\gT \models \varphi_{Q,0}(x_1,...,x_n)$.
    Recall that by lemma \ref{const} the value of $\gT \models
    \varphi_{Q,0}(x_1,...,x_n)$, is determined by $\langle Y; \leq^{\gT}, R_R, R_L \rangle$, the types
$\{tp_q^{MSO}((y,y'),\gT|_{T^3_{y,y'}}):y,y'\in Y,y<y',(y,y') \cap
Y = \emptyset\}$, and the types
$tp_q^{MSO}(c_{rt}^{\gT},\gT|_{T_L})$,
$tp_q^{MSO}(c_{rt}^{\gT},\gT|_{T_R})$ (see \ref{const} and
notation \ref{nota1}), but as $\gT$
    is correct these all are determined by $p$ so we are done.
    \end{proof}

    We can now conclude the proof of lemma \ref{tree const}. Define
$\gS = \{\bf{s}_u:u\subseteq\{1,...,k_2\}\}\cup\{\bf{s}^*\}$,
    and let $\gK$ be the constructible class of $\tau$-structures
defined by
    $\gP$ and $\gS$.
    Let $M$ be a $\tau$-structure in $\gK^{mo}_{\bf{c}}$. So we
    have $M \cong \gT_1^{[\bf{c}]}$ for some DB $k_2$-tree $\gT_1$. Let
$\gT_2$ be the
    correct extension of $\gT_1$.
    Let $A=\{x\in |M|: \gT_1 \models \varphi_{=,0}(x)\}$, and
    $\gA=\gT_2|_A\circledast_{\bf{s}^*}Null_{\emptyset}$
    From lemma \ref{ind} we have that $\gT_2^=$ is
    in the closer of $\gP$ under the operations in $\gS$ and hence
    so is $\gA$. From lemma \ref{s*} and the definition of
$\gT_1^{[\bf{c}]}$, we have
    that $\gA|_{\tau} = \gT_1^{[\bf{c}]} \cong M$, so $M \in \gK$ as
    desired.

\end{proof}

From lemmas \ref{tree const} and \ref{CONST IS TREE} we conclude
our main:

\begin{thm} \label{MAIN}

Let $\gK$ be a class of $\tau$-structures. Then $\gK$ is contained
in a $m$-constructible class for some $m\in\mathbb{N}$ iff $\gK$
is contained in a $3$-constructible class.

\end{thm}

The same holds for patch-width:

\begin{corollary}

Let $\tau$ be a nice vocabulary and $\gK$ a class of
$\tau$-structures. Then $\gK$ is contained in a class of bounded
$m$-ary patch-width for some $m\in\mathbb{N}$ iff $\gK$ is
contained in a class of bounded $3$-ary patch-width.

\end{corollary}

\begin{proof}
    Assume $\gK\subseteq\gK'$ for some $\gK'$ of bounded $m$-ary patch-width.
    By lemma \ref{PW IS CONST} $\gK'$ is $(m,0)$-constructible. By theorem \ref{MAIN}
    $\gK'$ is contained in some $3$-constructible $\gK''$. Notice
    that that the set $\gS$
    defined in the proof of \ref{tree const} satisfies that
    $\gS\subseteq\gS_{\tau^+,0,0,0}$ so $\gK''$ is in fact $(3,0)$-constructible.
    Notice further that in the proof of \ref{tree const} as we do not need null structures in the construction
    we may replace $\tau^+$ by a nice vocabulary.
    So by lemma \ref{CONST IS PW} $\gK''$ is a bounded $3$-ary
    patch-width class as desired.
\end{proof}

\section{A counter example for the unary case}
It turns out that we can not replace the number $3$ in theorem
\ref{MAIN} by $1$. This is because of the following:

\begin{thm} \label{EXAMPLE}
There exists a nice vocabulary $\tau$, and a class of
$\tau$-structures $\gK$, contained in some $3$-constructible
class, that is not contained in any $1$-constructible class.
\end{thm}

\begin{proof}
Let $\tau=\{R\}$ with $R$ 4-place relation symbol. Set $p \in
\mathbb{N}$ be large enough (to be defined later). Let $\gT$ be a
tree. For $x,y\in T$ Define:
\begin{itemize}

    \item
    $x \wedge^{\gT} y = x \wedge y =$ the $<^{\gT}$ minimal $z\in T$
    such that $z \geq^{\gT} x \text{ and } z \geq^{\gT} x$.

    \item
    $d^{\gT}(x,y)=d(x,y)=min \{|S|:S \subseteq T, x,x \wedge y \in
    S, \text{ S is dense in }(T,<^{\gT})\}$.

    \item
    $d^{\gT}_p(x,y)=d_p(x,y)=d(x,y) \pmod{p}$.

\end{itemize}

Let $q : \{0,...,p-1\}^2 \to \{0,1\}$ be some function that will
be defined later. We now define $\bf{c}$ a $0$-interpretation
scheme for $\tau$:
\begin{itemize}
    \item $k_1^{\bf{c}}=k_2^{\bf{c}}=0$.
    \item $\varphi_{=,0}^{\bf{c}}(x)=\lnot \exists y y>x$ i.e. the
    elements of the interpreted structure are the leafs of the tree.
    \item

$\varphi_{R,0}^{\bf{c}}(x_1,x_2,x_3,x_4)="q(d_p(x_1,x_2),d_p(x_3,x_4))=0"$.
\end{itemize}
We have to show that $\varphi_{R,0}$ is indeed a monadic formula
in $\tau_{trees}$. Note that there exists a monadic formula
$\varphi_{d_p=0}(x,y)$ such that for any tree $\gT$, $\gT \models
\varphi_{d_p=0}(x,y)$ iff $d^{\gT}_l(x,y)=0$.
$\varphi_{d_p=0}(x,y)$ will "say" that there exists a set $X$ such
that: \begin{itemize}
    \item $x,x \wedge y \in X$,
    \item if $z,z' \in X$ and $z<z''<z'$ then $z'' \in X$,
    \item if $z'$ is the immediate successor in $X$ of $z \in X$,
    then there exist exactly $p-1$ elements (of $T$) between them.
\end{itemize}
similarly we have formulas $\varphi_{d_p=i}(x,y)$ for $0 < i <p$.
Now define:
$$ \varphi_{R,0}(x_1,x_2,x_3,x_4)=\bigvee_{\substack{n_1,n_2 \in
\{0,...,p-1\} \\ n_1\equiv n_2 \pmod{p}}}
\varphi_{d_p=n_1}(x_1,x_2) \wedge \varphi_{d_p=n_2}(x_3,x_4).$$
This gives us $\bf{c}$ as desired. Define
$\gK=\gK^{mo,db}_{\bf{c}}$. By \ref{tree const} $\gK$ is contained
in a $3$-constructible class (in fact in a $3$-ary BPW class).
\newline For each $n \in \mathbb{N}$ let $M_n = ({}^{n
\geq}2,\triangleleft)$ i.e. $M_n$ is the complete binary tree of
depth $n$, and $N_n=M_n^{[\bf{c}]}$. Let $\gK'$ be a constructible
class of $\tau$-structures, so $\tau^+ = \tau_m$ for some $m \in
\mathbb{N}$. Towards contradiction assume that $N_n \in \gK'$ for
all $n \in \mathbb{N}$. Let $\gP$ be the set of "atomic"
structures associates with $\gK'$. w.l.o.g. we may assume that
$\gP$ consists of singleton structures only. Otherwise increase
$k$ by $max\{|M| : M \in \gP\}$ and construct each $M \in \gP$
from singletons of distinct colors. Now let $K \in \gK'$, and let
$(\gT,\gM)$ be a full representation of $K$ (see \ref{rep}).
Assume $K\cong N_n$ for some $n$. So we have a 1-1 function $f$,
from ${}^n2$ to the leafs of $\gT$, as every $\eta \in
{}^\frac{}{}n2$ corresponds to a unique element $a \in K$ under
the isomorphisms, and for every element of $a \in K$ there exist a
unique $t$ a leaf of $\gT$ such that $a = |M_t|$. Define
$f(eta)=t$. Note that $f$ is not onto as some of the leafs of
$\gT$ may be omitted during the creation process. For each $t \in
T$ let $A_t = \{f^{-1}(s) : s \leq^{\gT} \wedge s \in range(f)
\}$. So $A_t \subseteq {}^n2$. For each $eta \in A_t$ let
$a=a_{eta}=|M_{f(\eta)}|$. $a_{\eta}$ is an element of $M_t$, So
$A_t$ is divided into $2^k$ parts according to the color of
$a_{\eta}$ in $M_t$, (more formally according to the type
$tp_{qf}^{\tau^+ \setminus \tau}(a_{eta},M_t)$. So we have $B_T
\subseteq A_t$ such that $|B_t| \geq \frac{|A_t|}{2^k}$, and all
the elements of $f^{-1}(B_t)$ have the color. Now define:

$$C_t=\{d_p^{N_n}(\eta ,\eta \wedge \nu):\eta,\nu\in B_t\} \subseteq
\{0,...,p-1\}$$

We have $\frac{|A_t|}{2^k} \leq |B_t| \leq 2^{|C_t|}$ For the
right-hand inequality use induction on $|C_t|$. Hence we conclude
$$|A_t| \leq |B_t| \cdot 2^k \leq 2^{|C_t|+k}$$

Now note that if $C_t \neq \{0,...,p-1\}$, then $C_t \leq
n-\lfloor\frac{n}{p}\rfloor$ and hence $|A_t| \leq 2^{|C_t|+k}
\leq 2^{n-\lfloor\frac{n}{p}}\rfloor +k$. We now consider two
cases:

\begin{itemize}
\item[\textbf{Case 1}] There exist $s \in T$ with two immediate
successors $t_1,t_2 \in T$ such that: $|A_{t_1}|,|A_{t_2}| >
2^{n-\lfloor\frac{n}{p}}\rfloor +k$.
\end{itemize}

According to what we saw above we have $C_1=C_2=\{0,...,p-1\}$. So
for $l \in \{1,2\}$ we have $\langle (\rho_{t_l,i},\nu_{t_l,i}:i
\in \{0,...,p-1\} \rangle$ such that:
\begin{itemize}
    \item[$(\alpha)$] $\{\rho_{t_l,i},\nu_{t_l,i}:i
\in \{0,...,p-1\}\}$ all have the same color in $M_{t_l}$.
    \item[$(\beta)$] $d_p^{N_n}(\rho_{t_l,i},\nu_{t_l,i})=i$ for
    all $i<p$.
\end{itemize}
Denote by $m$ the number of quantifier free types of couples in
the vocabulary $\tau$ (actually in our case $m=2^{(2^4)}$). Note
that $m$ does not depend on $p$. So for each $l \in
\{1,2\}$$\{0,...,p-1\}$ is deviled into $m$ parts according to the
type $tp_{qf}((\rho_{t_l,i},\nu_{t_l,i}),M_{t_l})$. we claim that
we can (a priori) choose $p$ (large enough) and $q$ in such a way
that we can find: $i_1,i_2,j_1,j_2$ such that for each $l \in
\{1,2\}$: $\rho_{t_l,i_l},\nu_{t_l,i_l})$ and
$\rho_{t_l,j_l},\nu_{t_l,j_l})$ have the same quantifier free type
in vocabulary $\tau$ in $M_{t_l}$, and on the other hand:
$q(i_1,i_2) \neq q(j_1,j_2)$. This is of course a contradiction as
the quantifier free type of $\rho_{t_l,i_l},\nu_{t_l,i_l})$ and
$\rho_{t_l,j_l},\nu_{t_l,j_l})$ in vocabulary $\tau_k$ in
$M_{t_l}$ determines the value of
$R(\rho_{t_l,i_l},\nu_{t_l,i_l}),\rho_{t_l,j_l},\nu_{t_l,j_l})$ in
$M_s$ and hance in $M_{c_{rt}^\gT}$. But this value is true iff
$q(i_1,i_2)=0$ in contradiction with $q(i_1,i_2) \neq q(j_1,j_2)$.
Way can we choose $p$ and $q$ as desired? For a given $p$ the
number of functions from $\{0,...p-1\}^2$ to $\{0,1\}$ such that
we can not choose as above (i.e. functions that "respects" some
partition of $\{0,...,p-1\}$ into $m$ parts is the number of
partitions $m^p \cdot m^p$, time the number of functions that
"respect" that partition $2^{m \cdot m}$, or $2^{2p \log(m) +
m^2}$. The total number of functions is $2^{p^2}$. So if we choose
(a priori) $p$ such that $p^2 > 2p \log(m) + m^2$ we can choose a
function $q$ as desired. \newline

Assume now that the assumption of \textbf{Case 1} does not hold.
Assume also that we have chosen $n$ large enough such that
$2^{\lfloor \frac{n}{p} \rfloor -k} >4$. In this case we can find
$t_0,t_1,...,t_d \in T$ such that :
\begin{itemize}
    \item $d \geq 5$.
    \item $t_0=c_{rt}^{\gT}$.
    \item $t_d$ is a leaf of $\gT$.
    \item For $0 \leq i <d$, $t_{i+1}$ is an immediate successor in
$\gT$,
    of $t_{1}$.
    \item For $0 \leq i <d$, denote by $s_{i+1}$ the immediate
    successor of $t_i$ different from $t_{i+1}$, then
    $|A_{s_{i+1}}| \leq 2^{\lfloor \frac{n}{p} \rfloor -k}$.
\end{itemize}
Note that for any $0 < i \leq d$: $\bigcup_{0 < j \leq i}
{A_{s_{j}}}$ and $A_{t_i}$ is a partition of $A_{c_{rt}^{\gT}}$,
and that $|A_{c_{rt}^{\gT}}|=2^n$. So we can find $0 < i^* \leq d$
such that $|\bigcup_{0 < j \leq i^*} {A_{s_{j}}}|,|A_{t_{i^*}}| >
2^{\lfloor \frac{n}{p} \rfloor -k}$. We proceed similarly to
\textbf{Case 1}. As there we can find $\langle
(\rho_{t_{i^*},i},\nu_{t_{i^*},i} \in A_{t_{i^*}}:i \in
\{0,...,p-1\} \rangle$ that satisfy $(\alpha)$ and $(\beta)$
above, and the same for $\langle (\rho_i,\nu_i :i \in
\{0,...,p-1\} \rangle$ where $\rho_i,\nu_i \in \bigcup_{0 < j \leq
i^*} {A_{s_{j}}}$. Again let $m$ denote the number of quantifier
free types of couples in the vocabulary $\tau$. This time we want
to choose $p$ and $q$ in such a way that we can find: $i,j_1,j_2$
such that: $\rho_{t_{i^*},j_1},\nu_{t_{i^*},j_1})$ and
$\rho_{t_{i^*},j_2},\nu_{t_{i^*},j_2})$ have the same quantifier
free type in vocabulary $\tau$ in $M_{t_{i^*}}$, and on the other
hand: $q(i,j_1) \neq q(i,j_2)$. Again this is a contradiction as
the quantifier free type of
$\rho_{t_{i^*},j_l},\nu_{t_{i^*},j_l})$ for $l \in \{1,2\}$
determines the value of
$R(\rho_{t_{i^*},j_l},\nu_{t_{i^*},j_l}),\rho_i,\nu_i)$ in
$M_{c_{rt}^\gT}$. Again this value is true iff $q(i,j_l)=0$ in
contradiction with $q(i_1,i_2) \neq q(j_1,j_2)$. Way can we choose
$p$ and $q$ as desired? For a given $p$ the number of functions
from $\{0,...p-1\}^2$ to $\{0,1\}$ such that we can not choose as
above is the number of partitions $m^p$, times the number of
functions that "respect" that partition $2^{m\cdot p}$, or $2^{p
\log (m) + m \cdot p}$. So if we choose $p$ such that $p^2
> p \log(m) + m \cdot p$ we can choose a function $q$ as desired.
Note that the function we used for the second case will also work
for the first case so we can use one definition of $q$. In both
cases we get a contradiction and the proof is complete.

\end{proof}

\bibliographystyle{plain}
\bibliography{mon}

\end{document}